\documentclass[11pt]{article}
\usepackage[font=small,labelfont=bf]{caption}
\usepackage{latexsym,amsfonts,amsthm,amsmath,amssymb}
\usepackage{hyperref}
\usepackage{xcolor}
\usepackage{verbatim}
\usepackage{enumerate}
\usepackage{tabularx,graphicx}

\begin{document}

\begin{center}
{\LARGE On Euler's Solution of the simple Difference Equation} \\

\bigskip\bigskip

{Alexander Aycock, Johannes-Gutenberg University Mainz \\ Staudinger Weg 9, 55128 Mainz \\ \url{aaycock@students.uni-mainz.de}}

\bigskip


\end{center}

\medskip

\begin{abstract}\noindent 
In this note we will discuss Euler's solution of the simple difference equation that he gave in his paper{\it ``De serierum determinatione seu nova methodus inveniendi terminos generales serierum"} \cite{E189} (E189:``On the determination of series or a new method of finding the general terms of series") and also present a derivation for the values of the Riemann $\zeta$-function at positive integer numbers based on Euler's ideas.
\end{abstract}

\medskip

\section{Introduction}

In his paper{\it ``De serierum determinatione seu nova methodus inveniendi terminos generales serierum"} \cite{E189} (E189:``On the determination of series or a new method of finding the general terms of series"), Euler, amongst other difference equations, gave a general solution of the simple difference equation:

\begin{equation}
    \label{eq: Simple Difference Equation}
    f(x+1)-f(x)=g(x).
\end{equation}
He had found a solution to (\ref{eq: Simple Difference Equation}) in form of the Euler-Maclaurin summation formula before, e.g., in his paper{\it ``Inventio summae cuiusque seriei ex dato termino generali"} \cite{E47} (E47:``Finding of a sum of a series from the given general term"). But whereas the Euler-Maclaurin summation formula is a particular solution and leads to an asymptotic series  for most choices of $g(x)$, his solution offered in \cite{E189} is the complete solution to (\ref{eq: Simple Difference Equation}) and contains the Euler-Maclaurin summation formula as a special case.

Therefore, in this note we will present Euler's solution of (\ref{eq: Simple Difference Equation}) (see section \ref{sec: Euler's Solution of the Simple Difference Equation}), address a conceptual error in Euler's approach (see section \ref{sec: Discussion of Euler's Solution}) and we will show how to correct it (see section \ref{subsec: Correction of Euler's Formula}). Furthermore, we argue that Euler could have corrected his formula himself applying results that he discovered after he wrote \cite{E189} (see section \ref{subsec: Discussion}). Finally, we will present a derivation of the formula for the values of the Riemann $\zeta$-function at positive integer numbers based on the solution to the simple difference equation (see section \ref{sec: Application}).

\section{Euler's Solution of the Simple Difference Equation}
\label{sec: Euler's Solution of the Simple Difference Equation}

Euler's general idea was to transform (\ref{eq: Simple Difference Equation}) into a differential equation of infinite order with constant coefficients and apply the procedure he had formulated for the finite order case earlier in his paper{\it ``Methodus aequationes differentiales altiorum graduum integrandi ulterius promota"} \cite{E188} (E188:``The method to integrate differential equations of higher degrees expanded further").  In that paper he outlined the following procedure:

Given the differential equation:
\begin{equation}
    \label{eq: Inhomogeneous Finite}
    \left(a_0+a_1 \dfrac{d}{dx}+a_2 \dfrac{d^2}{dx^2} +\cdots+ a_n \dfrac{d^n}{dx^n} \right)f(x)= g(x),
\end{equation}
with complex coefficients $a_1, a_2, \cdots, a_n$, Euler told us to first find the zeros with their multiplicity of the following expression:

\begin{equation*}
    P(z)=a_0+a_1 z+a_2 z^2 +\cdots+ a_n z^n.
\end{equation*}
Next, assume $z=k$ is a solution of $P(z)=0$. Then, if $k$ is a simple zero\footnote{In this note, we will only need the case of simple zeros and hence will only state the corresponding formula. In \cite{E188}, Euler stated all cases from order 1 to 4 explicitly.} of $P(z)$, the solution of (\ref{eq: Inhomogeneous Finite}) is given by the sum of all functions of the form:

\begin{equation}
    \label{eq: Inhomogeneous Finite Solution}
    f_k(x)= \dfrac{e^{kx}}{P'(k)} \int e^{-kx}g(x)dx.
\end{equation}
Note that the indefinite integral introduces a constant of integration. 

Let's apply this to (\ref{eq: Simple Difference Equation}) by transforming it into a differential equation first. By Taylor's theorem we have:

\begin{equation*}
    f(x+1)=f(x)+\dfrac{d}{dx}f(x)+\dfrac{1}{2}\dfrac{d^2}{dx^2}f(x)+\dfrac{1}{3!}\dfrac{d^3}{dx^3}f(x)+\cdots
\end{equation*}
such that (\ref{eq: Simple Difference Equation}) can be rewritten as

\begin{equation}
    \label{eq: Differential Simple Difference}
    \left(\dfrac{d}{dx}+\dfrac{1}{2}\dfrac{d^2}{dx^2}+\dfrac{1}{3!}\dfrac{d^3}{dx^3}+\cdots\right)f(x) =g(x).
\end{equation}
Thus, according to Euler's approach we need to find all zeros and their multiplicities of the expression:

\begin{equation}
\label{eq: P(z) infinite}
    P(z) = \dfrac{z}{1!}+\dfrac{z^2}{2!}+\dfrac{z^3}{3!}+\dfrac{z^4}{4!}+\cdots = e^z-1.
\end{equation}
The general zero of this equation is $z=\log(1)$. But having established that the complex logarithm is a multivalued function in his work{\it ``De la controverse entre Mrs. Leibnitz et Bernoulli sur les logarithmes des nombres négatifs et imaginaires"} \cite{E168} (E168:``On the controversy between Leibnitz and Bernoulli on logarithms of negative and imaginary number"), Euler knew that (\ref{eq: P(z) infinite}) has infinitely many solutions, namely - aside from the trivial $z=0$ - the solutions are

\begin{equation*}
    \pm 2 \pi i, \pm 4 \pi i, \pm 6 \pi i, \pm 8 \pi i, \cdots.
\end{equation*}
Therefore, the formula (\ref{eq: Inhomogeneous Finite Solution}) applied to (\ref{eq: Differential Simple Difference}) and hence the solution to (\ref{eq: Simple Difference Equation}) gives:

\begin{equation}
    \label{eq: Euler Solution Simple Difference}
    f(x) = \int g(x)dx + e^{-2 \pi ix}\int g(x) e^{2 \pi ix}dx +  e^{2 \pi ix}\int g(x) e^{-2 \pi ix}dx
\end{equation}
\begin{equation*}
   +  e^{-4 \pi ix}\int g(x) e^{4 \pi ix}dx +  e^{4 \pi ix}\int g(x) e^{-4 \pi ix}dx+ \cdots
\end{equation*}
This is the solution Euler gave in \cite{E189}. Unfortunately, it is not quite correct. We will discuss this in the following section.

\section{Discussion of Euler's Solution}
\label{sec: Discussion of Euler's Solution}

\subsection{Example of linearly increasing Differences}

Applying Euler's formula (\ref{eq: Euler Solution Simple Difference}) to certain examples, we quickly discover that it does not give the correct results. For the purpose of illustration, let us take $g(x)=x$ such that we want to solve:

\begin{equation}
\label{eq: Example Linear Difference}
    f(x+1)-f(x)=x.
\end{equation}
The general solution to this equation is easily seen to be given as

\begin{equation}
\label{eq: Solution Linear Difference}
    f(x)= \dfrac{1}{2}x(x-1)+ h(x),
\end{equation}
where $h(x)$ satisfies $h(x+1)=h(x)$. Now let us apply (\ref{eq: Euler Solution Simple Difference}). 
For this, we need to evaluate:

\begin{equation*}
    e^{-2k\pi ix}\int x e^{2k\pi i x}dx=  \dfrac{1-2k \pi ix}{4 \pi^2 k^2}+C_k  e^{-2k\pi ix},
\end{equation*}
where $C_k$ is a constant of integration. For $k=0$, we have $\int xdx=\frac{x^2}{2}+C_0$, where $C_0$ is the constant of integration. Inserting all this into (\ref{eq: Euler Solution Simple Difference}), we find:

\begin{equation*}
    f(x) =\dfrac{x^2}{2}+C_0 + \sum_{k\in \mathbb{Z}\setminus \lbrace 0 \rbrace} \left(\dfrac{1}{4 \pi^2 k^2}-\dfrac{x}{2k \pi i}\right)+ C_k  e^{-2k\pi ix}.
\end{equation*}
Calling $C_0+\sum_{k\in \mathbb{Z}\setminus \lbrace 0 \rbrace}  C_k  e^{-2k\pi ix}=h(x)$, we see that $h(x)=h(x+1)$. Furthermore,

\begin{equation*}
     \sum_{k\in \mathbb{Z}\setminus \lbrace 0 \rbrace} \dfrac{1}{2k \pi i} =0,
\end{equation*}
because all terms cancel. Finally,

\begin{equation*}
    \sum_{k\in \mathbb{Z}\setminus \lbrace 0 \rbrace} \dfrac{1}{4 \pi^2 k^2} = \dfrac{2}{4 \pi^2}\sum_{k=1}^{\infty} \dfrac{1}{k^2}= \dfrac{1}{2\pi^2}\cdot \dfrac{\pi^2}{6}= \dfrac{1}{12},
\end{equation*}
where we used the result $\sum_{k=1}^{\infty}\frac{1}{k^2}=\frac{\pi^2}{6}$ that Euler had discovered in his paper{\it ``De summis serierum reciprocarum"} \cite{E41} (E41:``On the sums of series of reciprocals") in the last step.

Thus, Euler's formula (\ref{eq: Euler Solution Simple Difference}) gives the following solution to (\ref{eq: Example Linear Difference}):

\begin{equation*}
    g(x) = \dfrac{x^2}{2}+h(x),
\end{equation*}
where we absorbed the value $\frac{1}{12}$ in the formula in the periodic function. Comparing this result to (\ref{eq: Solution Linear Difference}), we see that the solution from Euler's formula is off by the term $-\frac{x}{2}$. In the next sections we will elaborate on why (\ref{eq: Euler Solution Simple Difference}) is wrong and how to correct it.

\subsection{Correction of Euler's Formula}
\label{subsec: Correction of Euler's Formula}

Euler's formula (\ref{eq: Euler Solution Simple Difference}) is actually almost correct. Indeed, the correct formula reads:

\begin{equation}
    \label{eq: Simple Difference Correct Solution}
    f(x) ={\color{red}-\dfrac{1}{2}g(x)}+ \int g(x)dx + e^{-2 \pi ix}\int g(x) e^{2 \pi ix}dx +  e^{2 \pi ix}\int g(x) e^{-2 \pi ix}dx
\end{equation}
\begin{equation*}
     e^{-4k \pi ix}\int g(x) e^{4 \pi ix}dx +  e^{4 \pi ix}\int g(x) e^{-4 \pi ix}dx+ \cdots
\end{equation*}
such that Euler's formula is off by just the term $-\frac{1}{2}g(x)$. As we mentioned in the introduction, Euler missed this term since the method of construction the solution to a differential equation from the zeros of the characteristic polynomial does not carry over smoothly from the finite to the infinite order case. Indeed, we have to construct the solution from the reciprocal of the characteristic polynomial, if we want the method to be applicable in the infinite order case. For, setting $z=\frac{d}{dx}$, we can rewrite (\ref{eq: Differential Simple Difference}) as:

\begin{equation}
\label{eq: Solution formal}
    f(x) = \dfrac{1}{P(z)}g(x)
\end{equation}
with $P(z)=e^z-1$. In order to apply the operator $\frac{1}{P(z)}$ to $g(x)$, we need to rewrite it in integer powers of $z$. There are many ways to achieve this task. The one we will need to prove (\ref{eq: Simple Difference Correct Solution}) is the following partial fraction decomposition that can be proved, e.g., by using complex analysis\footnote{In section \ref{subsec: Discussion} we will present a proof that uses only method that were available to Euler.}:

\begin{equation}
    \label{eq: Partial Fraction Decomposition}
    \dfrac{1}{e^z-1}=-\frac{1}{2}+ \sum_{k \in \mathbb{Z}\setminus \lbrace 0\rbrace}^{\infty} \dfrac{1}{z-2k\pi i}.
\end{equation}
Thus, next we have to evaluate:

\begin{equation}
    \label{eq: Single Term}
    \dfrac{1}{z-2 k \pi i}g(x).
\end{equation}
Writing $2 k \pi i = \alpha$ we have:

\begin{equation*}
    \dfrac{1}{z-\alpha}g(x) = \dfrac{1}{z \left(1- \frac{\alpha}{z}\right)}g(x) = \sum_{n=0}^{\infty}\dfrac{\alpha^n}{z^{n+1}}g(x).
\end{equation*}
Since $z=\dfrac{d}{dx}$, we can interpret $\frac{1}{z}$ as an integral and hence $\dfrac{1}{z^n}$ is an $n$-times iterated integral. Writing $\int^n$ for a $n$-times iterated integral, the following formula holds:

\begin{equation}
    \label{eq: Iterated integral}
    \int^n g(x)dx = \int\limits_{}^{x}\dfrac{(x-t)^{n-1}}{(n-1)!}g(t)dt.
\end{equation}
Inserting this into (\ref{eq: Single Term}), we have:

\begin{equation*}
     \dfrac{1}{z-\alpha}g(x)= \sum_{n=0}^{\infty} \alpha^n \int\limits_{}^{x}\dfrac{(x-t)^{n}}{n!}g(t)dt = \int\limits_{}^{x} e^{\alpha (x-t)}g(t)dt.
\end{equation*}
Therefore, by (\ref{eq: Partial Fraction Decomposition}) our equation (\ref{eq: Solution formal}) reads:

\begin{equation*}
f(x) = -\dfrac{1}{2}g(x)+\sum_{k\in \mathbb{Z}} \int\limits_{}^{x} e^{2 k\pi i (x-t)}g(t)dt =  -\dfrac{1}{2}g(x)+\sum_{k\in \mathbb{Z}}e^{2 k \pi i x} \int\limits_{}^{x} e^{-2 k\pi i t}g(t)dt,
\end{equation*}
which is (\ref{eq: Simple Difference Correct Solution}). It is the same solution as in \cite{We14}, which derived (\ref{eq: Simple Difference Correct Solution}) using complex analysis.

\subsection{Discussion}
\label{subsec: Discussion}

Although we operated on a purely formal basis in our derivation of (\ref{eq: Simple Difference Correct Solution}), the procedure can be justified applying the Fourier transform which allows to consider (\ref{eq: Differential Simple Difference}) as an algebraic equation in the new variable, say, $p$. To find (\ref{eq: Single Term}) we then need the inverse Fourier transform, which we can either calculate using complex analysis or look up in a table.

But Fourier analysis was not available to Euler, of course. Nevertheless, we argue that Euler could have given the proof we presented himself. The proof hinges essentially on the proof of (\ref{eq: Partial Fraction Decomposition}). Later in his career, in his paper{\it ``De resolutione fractionum transcendentium in infinitas fractiones simplices"} \cite{E592} (E592:``On the resolution of transcendental fractions into infinitely many simple fractions"), Euler indeed considered partial fraction decompositions of transcendental functions. The method outlined there would have given him the formula:

\begin{equation}
\label{eq: Partial Fraction Euler}
    \dfrac{1}{e^z-1}= R(z)+\sum_{k \in \mathbb{Z}} \dfrac{1}{z-2 k \pi i},
\end{equation}
where $R(z)$ is a function to be determined. Next, one could expand the sum into a Laurent series around $z=0$ by expanding each geometric series and compare it to the Laurent series obtained by direct expansion. The direct expansion reads:

\begin{equation}
    \label{eq: Direct Expansion}
    \dfrac{1}{e^z-1}=-\dfrac{1}{2}+\dfrac{1}{z}+\sum_{k=0}^{\infty} B_n\dfrac{z^n}{n!},
\end{equation}
where $B_n$ are the Bernoulli numbers. Since Euler considered a similar function, namely $\frac{z}{1-e^{-z}}$, and its series expansion around $z=0$ in his work{\it ``De seriebus quibusdam considerationes"} \cite{E130} (E130:``Considerations about certain series"), the previous formula could definitely also been found by him. Finally, comparing the Laurent series obtained from (\ref{eq: Partial Fraction Euler}) to (\ref{eq: Direct Expansion}), we can infer that $R(z)=-\frac{1}{2}$.


\section{An Application of the Solution to the simple Difference Equation}
\label{sec: Application}

In this section, we want to consider the choice $g(x)=x^n$ for $n \in \mathbb{N}$ in (\ref{eq: Simple Difference Equation}), since it is one of the few cases in which (\ref{eq: Simple Difference Correct Solution}) can be evaluated explicitly. As it will turn out, we will be led to the values $\zeta(2n)$, i.e., the sums

\begin{equation}
\label{eq: zeta(2n)}
    \zeta(2n):=\sum_{k=1}^{\infty} \dfrac{1}{k^{2n}}
\end{equation}
in the process. Euler evaluated these sums on many occasions using a large number of different methods. We mention his papers \cite{E41} and \cite{E130} as examples, but the way we will arrive at those values seems to be different from all methods used by Euler.

\subsection{Preparation}

Considering (\ref{eq: Simple Difference Correct Solution}) we need to evaluate the expression: $ e^{a x} \int e^{-ax}x^n dx$. This can be done as follows: First, we note that

\begin{equation}
\label{eq: Case x^0}
    \int e^{-ax}dx = -\dfrac{e^{ax}}{a}=- e^{ax}\cdot a^{-1},
\end{equation}
where we omitted the constant of integration, since it will not be necessary in the following. Next, we differentiate (\ref{eq: Case x^0}) with respect to $a$ exactly $n$ times. The left-hand side gives:

\begin{equation*}
    \dfrac{d^n}{da^n} \int e^{-ax}dx= \int  \dfrac{d^n}{da^n} e^{-ax}dx  =(-1)^n \int e^{-ax}x^ndx,
\end{equation*}
whereas the right-hand side gives:

\begin{equation*}
   - \dfrac{d^n}{da^n}  e^{-ax}\cdot a^{-1} = - \sum_{k=0}^n \binom{n}{k} \dfrac{d^k}{da^k} e^{-ax}\cdot \dfrac{d^{n-k}}{da^{n-k}} a^{-1}
\end{equation*}
\begin{equation*}
    = -e^{-ax}\cdot \dfrac{(-1)^n}{a^{n+1}}\sum_{k=0}^n \dfrac{n!}{k!}(ax)^k,
\end{equation*}
where we used Leibniz' rule for the differentiation of products in the first step. Thus, combining both results we arrive at:
\begin{equation*}
    e^{a x} \int e^{-ax}x^n dx = -\dfrac{1}{a^{n+1}}\sum_{k=0}^{n} \dfrac{n!}{k!}a^kx^k.
\end{equation*}
Inserting this into (\ref{eq: Simple Difference Correct Solution}) for the special case $g(x)=x^n$ we get:

\begin{equation}
    \label{eq: Solution x^n}
    f(x)= \dfrac{x^{n+1}}{n+1}-\dfrac{x^n}{2}- \sum_{k \in \mathbb{Z}\setminus \lbrace 0 \rbrace} \dfrac{1}{(2k \pi i)^{n+1}}\cdot \sum_{j=0}^n \dfrac{n!}{j!}(2k\pi i)^j x^j+h(x),
\end{equation}
where $h(x)$ satisfies $h(x+1)=h(x)$.

\subsection{The Application}

(\ref{eq: Solution x^n}) is the general solution to (\ref{eq: Simple Difference Equation}) for the particular choice $g(x)=x^n$. But we can also easily find a particular solution to (\ref{eq: Simple Difference Equation}) by noting that for integer $x$:

\begin{equation*}
    f(x)= \sum_{k=1}^{x-1} g(k) = \sum_{k=1}^{x}g(k) - g(x)
\end{equation*}
satisfies the equation. Therefore, for the particular choice $g(x)=x^n$ we also have the solution:

\begin{equation}
\label{eq: Faulhaber}
    f(x) = \sum_{k=1}^{x-1}k^n= \sum_{k=1}^{x}k^n - x^n = \dfrac{x^{n+1}}{n+1}+\dfrac{x^n}{2}+\dfrac{1}{n+1}\sum_{j=2}^{n} \binom{n+1}{j}B_j x^{n+1-j}-x^n
\end{equation}
\begin{equation*}
    = \dfrac{x^{n+1}}{n+1}-\dfrac{x^n}{2}+\dfrac{1}{n+1}\sum_{j=2}^{n} \binom{n+1}{j}B_j x^{n+1-j},
\end{equation*}
where we used Faulhaber's formula for the sums of integer powers and $B_n$ is the $n$-th Bernoulli number as above.  (\ref{eq: Faulhaber}) is a polynomial in $x$ and hence $x$ is not restricted to integer values in this form. Let us transform (\ref{eq: Solution x^n}) into a similar form. Ignoring the periodic function we have:

\begin{equation}
\label{eq: Ordered x^n}
    f(x)= \dfrac{x^{n+1}}{n+1}-\dfrac{x^n}{2}-  \sum_{j=0}^n \dfrac{n!}{j!} \sum_{k \in \mathbb{Z}\setminus \lbrace 0 \rbrace}(2k\pi i)^{j-(n+1)} x^j.
\end{equation}
Since (\ref{eq: Faulhaber}) and (\ref{eq: Ordered x^n}) differ only by a periodic function, we can compare coefficients of respective powers of $x$. Let us call

\begin{equation}
    \label{eq: Definition Coefficients B}
    B(n,j)= \dfrac{1}{n+1}\binom{n+1}{j}B_j \quad \text{for} \quad j \geq 2,
\end{equation}
for all other values of $j$ we set $B(j,n)=0$; furthermore, we set

\begin{equation}
    \label{eq: Definition Coefficients A}
A(n,j)= - \dfrac{n!}{j!}\sum_{k \in \mathbb{Z}\setminus \lbrace 0 \rbrace}(2k\pi i)^{j-(n+1)}.
\end{equation}
Then, comparing coefficients from (\ref{eq: Faulhaber}) and (\ref{eq: Ordered x^n}) gives:

\begin{equation*}
     A(n,n+1-j)=B(n,j).
\end{equation*}
Thus, substituting the values form (\ref{eq: Definition Coefficients A}) and (\ref{eq: Definition Coefficients B}), respectively:

\begin{equation*}
     - \dfrac{n!}{(n+1-j)!}\sum_{k \in \mathbb{Z}\setminus \lbrace 0 \rbrace}(2k\pi i)^{-j}= \dfrac{1}{n+1}\binom{n+1}{j}B_j.
\end{equation*}
Finally, solving for the sum:

\begin{equation*}
    \sum_{k \in \mathbb{Z}\setminus \lbrace 0 \rbrace}(2k\pi i)^{-j}=- \dfrac{(n+1-j)!}{n!}\cdot \dfrac{1}{n+1}\cdot \dfrac{(n+1)!}{j!(n+1-j)!}B_j=-\dfrac{B_j}{j!}.
\end{equation*}
Thus,

\begin{equation*}
     \sum_{k \in \mathbb{Z}\setminus \lbrace 0 \rbrace}k^{-j} = -(2\pi i)^{j}\dfrac{B_j}{j!}.
\end{equation*}
But due canceling terms, the sum vanishes for odd $j$ such that we arrive at:

\begin{equation}
\label{eq: Formula zeta(2n)}
    \zeta(2j)=\sum_{k=1}^{\infty} \dfrac{1}{k^{2j}}=\dfrac{(-1)^{j-1}(2\pi)^{2j} B_{2j}}{2(2j)!}.
\end{equation}
This is Euler's famous formula for the even values of the $\zeta$-function that he gave, e.g., in \cite{E130}.

\section{Conclusion}

In this note we considered Euler's general solution to the simple difference equation (\ref{eq: Simple Difference Equation}) that he gave in \cite{E189}. His final formula (\ref{eq: Euler Solution Simple Difference}) is slightly incorrect due to unjustified application of his solution  (\ref{eq: Inhomogeneous Finite Solution}) to differential equations of infinite order. Nevertheless, we discussed how to fix Euler's derivation (see section \ref{subsec: Correction of Euler's Formula}) and also argued that Euler could have done so himself, if he just reconsidered the same subject later in his career (see section \ref{subsec: Discussion}).  Furthermore, we used the correct solution (\ref{eq: Simple Difference Correct Solution}) to (\ref{eq: Simple Difference Equation}) to give a proof of Euler's famous formula for the values of the Riemann $\zeta$-function at even positive integers (see section \ref{sec: Application}). The method of derivation of (\ref{eq: Formula zeta(2n)}) that we presented seems to be not to have been used by Euler in any of his other papers. Finally, we mention that our approach allowed to derive the exact values of (\ref{eq: zeta(2n)}) from the corresponding {\it finite} sums of natural powers (\ref{eq: Faulhaber}). Although this is clear, since (\ref{eq: Faulhaber}) and (\ref{eq: Formula zeta(2n)}) are connected via the Bernoulli numbers -- as Euler also pointed out, e.g., in \cite{E130} --, the deeper explanation for this connection is provided by (\ref{eq: Simple Difference Correct Solution}).

Despite the minor mistake, \cite{E189} is an interesting paper and contains subjects and approaches that are not found in any other of Euler's papers.


\end{document}